\documentclass[12pt]{article}
\usepackage{theorem,amsfonts}
\newtheorem{theorem}{Theorem}[section]
\newtheorem{proposition}{Proposition}[section]
\newtheorem{lemma}{Lemma}[section]
\newtheorem{corollary}{Corollary}[section]
\theorembodyfont{\upshape}
\newtheorem{definition}{Definition}
\newtheorem{remark}{Remark}

\newtheorem{example}{Example}[section]
\newtheorem{proof}{Proof}

\newcommand{\bt}{\begin{theorem}}
\newcommand{\et}{\end{theorem}}
\newcommand{\bl}{\begin{lemma}}
\newcommand{\el}{\end{lemma}}
\newcommand{\bp}{\begin{proposition}}
\newcommand{\ep}{\end{proposition}}
\newcommand{\bex}{\begin{example}}
\newcommand{\eex}{\end{example}}
\newcommand{\bc}{\begin{corollary}}
\newcommand{\ec}{\end{corollary}}
\newcommand{\bo}{\begin{proof}}
\newcommand{\eo}{\end{proof}}
\newcommand{\bd}{\begin{definition}}
\newcommand{\ed}{\end{definition}}
\newcommand{\br}{\begin{remark}}
\newcommand{\er}{\end{remark}}
\newcommand{\be}{\begin{enumerate}}
\newcommand{\ee}{\end{enumerate}}

\begin{document}

\def\frontpage#1#2#3#4#5{\begin{titlepage}
 \let\footnotesize\small      
 \let\footnoterule\relax      
 \setcounter{page}{0}
   {{}\hfill #1 \\ \mbox{}\hfill  #2 \\ \mbox{}\hfill  
  \texttt{http://www.isibang.ac.in/\raisebox{-.7ex}{$\widetilde{\hspace{.5em}}$}\,statmath/eprints}} 
  \vfill
 \begin{minipage}{450pt}                 
   \begin{center}
     {\huge #3}
        \vskip 2em 
    {\large                  
     \begin{tabular}[t]{c}#4     
     \end{tabular}\par}
 \end{center} 
 \end{minipage}
   \vfill
\begin{minipage}{450pt}
\begin{center}{\Large #5}
\end{center} \par
\end{minipage}
\end{titlepage}
\setcounter{footnote}{0}       
\let\frontpage\relax}


\frontpage{isibang/ms/2009/2}{April 6th, 2009}{On the existence of ergodic 
automorphisms in ergodic ${\mathbb Z} ^d$-actions on compact groups}
{\sc C. R. E. Raja}
{{Indian Statistical Institute, Bangalore Centre}\\ 
{8th Mile Mysore Road, Bangalore, 560059 India}}
\newpage
\mbox{}
\thispagestyle{empty}
\newpage
\setcounter{page}{1}

\title{On the existence of ergodic automorphisms in ergodic 
${\mathbb Z} ^d$-actions on compact groups}
\author{C. R. E. Raja}
\date{}
\maketitle

\let\ol=\overline 
\let\epsi=\epsilon 
\let\vepsi=\varepsilon 
\let\lam=\lambda 
\let\Lam=\Lambda 
\let\ap=\alpha 
\let\te=\tilde 
\let\vp=\varphi 
\let\ra=\rightarrow 
\let\Ra=\Rightarrow 
\let \Llra=\Longleftrightarrow 
\let\Lla=\Longleftarrow 
\let\lra=\longrightarrow 
\let\Lra=\Longrightarrow 
\let\ba=\beta 
\let\ov=\overline 
\let\ga=\gamma 
\let\Ba=\Delta 
\let\Ga=\Gamma 
\let\Da=\Delta 
\let\Oa=\Omega 
\let\Lam=\Lambda 
\let\un=\upsilon

\newcommand{\bR}{{\bf R}}
\newcommand{\bS}{{\bf S}}
\newcommand{\bn}{{\bf n}}
\newcommand{\cM}{{\cal M}}
\newcommand{\cN}{{\cal N}}
\newcommand{\Z}{{\mathbb Z}}
\newcommand{\Q}{{\mathbb Q}}
\newcommand{\N}{{\mathbb N}}
\newcommand{\R}{{\mathbb R}}
\newcommand{\C}{{\mathbb C}}
\newcommand{\T}{{\mathbb T}}
\newcommand{\F}{{\mathbb F}}

\begin{abstract}
Let $K$ be a compact metrizable group and $\Ga$ be a finitely generated group of 
commuting automorphisms of $K$.  We show that ergodicity of $\Ga$ implies $\Ga$ 
contains ergodic automorphisms if center of the action, 
$Z(\Ga ) = \{ \ap \in {\rm Aut}~(K) \mid \ap {\rm~~ commutes ~~with ~~ 
elements ~~ of ~~} \Ga \}$ has DCC.  To explain that the condition on the 
center of the action is not restrictive, we discuss certain abelian groups 
which in particular, retrieves Theorems of Berend \cite{Be} and 
Schmidt \cite{Sc1} proved in this context.
\end{abstract} 



\begin{section}{Introduction}

We shall be considering actions on compact groups.  All compact groups considered 
in this article are assumed to be metrizable and all automorphisms are assumed 
to be continuous.  Let $K$ be a compact group and ${\rm Aut }~(K)$ be the group 
of all continuous automorphisms of $K$.  

\bd
Let $K$ be a compact group and $\omega _K$ be the normalized Haar measure on $K$.  
Let $\Ga$ be a group of automorphisms of $K$.  We say that $\Ga$ is ergodic on 
$K$ if any $\Ga$-invariant Borel set $A$ of $K$ satisfies $\omega _K(A) =0$ or 
$\omega _K(A) =1$.  We say that $\ap \in {\rm Aut}~(K)$ is ergodic on $K$ if the 
action of $\{ \ap ^n \mid n\in \Z \}$ on $K$ is ergodic.
\ed

Ergodic action has been characterized in many forms, we now recall two such 
characterizations due to Berend (cf. Theorem 2.1 of \cite{Be}):  $\Ga$ is ergodic 
on $K$ if and only if there exists $x \in K$ such that the orbit $\Ga (x)$ is 
dense in $K$ if and only if for any non-trivial irreducible unitary representation 
$\pi$, the orbit $\{ \ap (\pi ) \mid \ap \in \Ga \}$ has infinitely many 
inequivalent unitary representations where $\ap (\pi )$ is given by 
$\ap (\pi )(x) = \pi (\ap ^{-1}(x))$ for all $x \in K$.  Also, 
Berend \cite{Be} proved the following result.  

\vskip 0.1in

\noindent {\bf Berend's Theorem:}  Ergodic action of a group $\Ga$ of 
commuting epimorphisms on a compact connected finite-dimensional abelian group 
contains ergodic automorphisms.

\vskip 0.1in

Recently Berend's Theorem was proved for certain hereditarily ergodic actions of 
solvable groups on compact connected finite-dimensional abelian groups (cf. 
\cite{BG}) and for ergodic actions of nilpotent groups on compact connected 
finite-dimensional abelian groups (cf. \cite{Ra}): this type of result was called 
a local to global 
correspondence in \cite{Ra} as no $\ap$ in a group $\Ga$ is ergodic implies the 
whole group $\Ga$ is not ergodic.  No analogue of Berend's Theorem is known for 
actions on general compact groups.  We will now look at the following examples in 
this context.  

\bex\label{e1}
Let $M$ be a compact group and $\Ga$ be a countably infinite group.  Then the 
(left) shift action of $\Ga$ on $M^\Ga$ is defined for $\ap \in \Ga$ and $f 
\in M^\Ga$, by $$\ap f(\ba) = f(\ap ^{-1}\ba)$$ for all $\ba \in \Ga$.  Then it can 
be seen that the shift action of $\Ga$ on $M^\Ga$ is ergodic (cf. \cite{Ra} 
when $M$ is abelian).  
If $\Ga$ is a torsion group (for example, we may take $\Ga$ to be the group of 
all finite permutations), then no element of $\Ga$ is ergodic.  Thus, there are 
ergodic actions without having any ergodic automorphisms.
\eex

If $\Ga$ is assumed to be a finitely generated solvable group, then $\Ga $ is 
torsion implies $\Ga$ is finite.  Hence actions as in Example \ref{e1} does not 
exist if $\Ga$ is assumed to be a finitely generated solvable group.  However the 
following example shows that even if we assume $\Ga$ to be finitely generated 
abelian, Berend's Theorem need not be true for any compact group.

\bex\label{e2} Let $M$ be a compact group with an ergodic automorphisms $\tau$.  
For any $(i, j) \in \Z ^2 \setminus (0, 0)$, define $K_{i, j} = M$ with $\Z 
^2$-action defined by $(n, m)\mapsto \tau ^{mi-nj}$ and define $K_{0,0} = 
M^{\Z^2}$ with shift action of $\Z ^2$.  Let $K= \prod _{i, j} K_{i, j}$ and 
define coordinate-wise $\Z ^2$-action on $K$.  Then it can be verified that this 
action is faithful, that is no non-zero $(n ,m)$ acts trivially on $K$.  We first 
claim that $\Z ^2$-action on $K$ is ergodic.  It is sufficient to claim that $\Z 
^2$-action on each $K_{i, j}$ is ergodic.  Let $(i, j)\in \Z ^2 \setminus (0,0)$.  
Since $\tau$ is ergodic on $M$, if $j \not = 0$, then $(i+1, j)$ whose action 
on $K_{i, j}$ is by $\tau ^{-j}$ is ergodic on $K_{i,j}$ and 
if $i\not =0$, then $(i, j+1)$ whose action on $K_{i,j}$ is by $\tau ^i$ is 
ergodic on $K_{i,j}$.  We know that the shift action of 
$\Z ^2$ on $K_{0,0} = M^{\Z ^2}$ is ergodic.  Thus, $\Z ^2$-action on $K$ is 
ergodic.

We next claim that no element $(i,j)$ of $\Z ^2$ is ergodic on $K$.  It is clear 
that $(0,0)$ is not ergodic on $K$.  Now take any $(i, j) \in \Z ^2 \setminus 
(0,0)$.  Now, $(i,j)$-action on $K_{i,j}$ is given by $\tau ^{ji-ij}$, identity.  
That is, $(i, j)$ is not ergodic on $K_{i, j}$.  Since $K_{i, j}$ is a factor of 
$K$, we get that $(i, j)$ is not ergodic on $K$.  Thus, no $(n,m) \in \Z^2$ 
is ergodic on $K$.  
\eex  

In the above example the $\Z ^2$-action on the compact group has infinitely many 
invariant normal subgroups resulting in infinite chain of closed normal 
invariant subgroups.  Motivated by example \ref{e2} we are lead to consider the 
following type of actions. 

\bd 
Action of a group $\Ga$ on a compact group $K$ is said to have descending chain 
condition (abbr. as DCC) if any decreasing sequence $(K_i)_{i\geq 1}$ of 
$\Ga$-invariant closed subgroups is finite, that is, there exists a $n\geq 1$ 
such that $K_n = K_{n+i}$ for all $i$.  
\ed

The condition DCC was studied in the context of ergodic actions and other type of 
group actions (cf. \cite{KS}).  Kitchens and Schmidt \cite{KS} proved that 
action of countable abelian groups is a projective limit of actions having DCC 
and action of a finitely generated abelian group $\Ga$ on a compact group $K$ has 
DCC if and only if $K$ is isomorphic to a shift-invariant subgroup of $M^\Ga$ for 
some compact real Lie group $M$ (cf. Theorems 3.2 and 3.16 of \cite{KS}): see 
Example \ref{e1} for definition of the shift action of $\Ga$ on $M^\Ga$.  It can 
be easily seen that $K$ and $\Ga$ in Example 
\ref{e2} does not satisfy DCC: take $K_n = \prod _{i+j\geq n} K_{i+j}$, $(K_n)$ 
is a decreasing sequence of closed (normal) $\Ga$-invariant subgroups 
but $K_n \not = K_m$ for all $n \not = m$.  

Schmidt \cite{Sc1} proved that if $\Ga$ is a finitely generated abelian 
group of automorphisms of a compact abelian group $K$ such that $\Ga$ has DCC on $K$, 
then $\Ga$ is ergodic on $K$ if and only if $\Ga$ contains ergodic automorphisms 
(cf. Theorem 3.3 and Corollary 3.4 of \cite{Sc1}).  
 
We now introduce center of actions.  If a group $\Ga$ acts on a compact group 
$K$, then center of the action denoted by $Z(\Ga )$ is defined to be $Z(\Ga ) 
= \{ \ap \in {\rm Aut}~(K) \mid \ap \ba = \ba \ap ~~{\rm for ~~all}~~\ba \in \Ga 
\}$.  If $\Ga$ is generated by $\ap _1, \cdots , \ap _n \in \Ga$ we may write 
$Z(\ap _1, \ap _2 , \cdots ,\ap _n)$ instead of $Z(\Ga )$.

In this article we study ergodic actions on compact groups and prove that ergodic 
actions of finitely generated abelian groups on compact groups contain ergodic 
automorphisms,   
hence extending Berend's Theorem provided that center of the action has DCC.  
Since a finitely generated abelian 
group $\Ga$ has a subgroup $\Ga '$ of finite index isomorphic to $\Z ^d$ for some 
$d\geq 0$ and $\Ga$ is ergodic if and only if any subgroup of finite index is  
ergodic, our results are about ergodic $\Z^d$-actions.   

We provide examples of compact groups where the center of certain type of 
automorphisms have DCC, 
as a consequence ergodic actions of finitely generated abelian groups of 
automorphisms of such type contain ergodic automorphisms.  This in particular, 
retrieves Berend's Theorem quoted above and Theorem 3.3 of Schmidt 
\cite{Sc1} which explains that the condition on center of the action is not 
restrictive.  

We will now see certain type of action which is needed for our purpose.  

\bd
We say that the action of a group $\Ga$ on a compact group $K$ is distal if 
for any $x\in K \setminus (e)$, $e$ is not in the closure of the orbit 
$\Ga (x) = \{ \ap (x) \mid \ap \in \Ga \}$.  We say that $\ap \in 
{\rm Aut}~(K)$ is distal on $K$ if the action of $\{ \ap ^n \mid n\in \Z \}$ 
on $K$ is distal.

\ed

Distal actions were studied by many in different contexts (cf. \cite{El}, 
\cite{Fu} and \cite{Ra}).  As in \cite{Ra} distality plays a useful role in our 
proofs.  It is a well-known interesting fact that if $\Ga$ is distal on $K$, then 
$\Ga$ is distal on $K/L$ for any $\Ga$-invariant normal subgroup $L$ of $K$ (cf. 
\cite{El}) and we quote this result as we use this result often without 
quoting.  
\end{section}

\begin{section}{Ergodic actions}

In this section we prove a few basic results on ergodic actions that are 
needed for our purpose and we provide a 
proof of these results as we could not locate these results in the literature.  

Let $K$ be a compact group.  If $\pi$ is a continuous unitary representation of 
$K$, then for any automorphism $\ap$ of $K$, 
$\ap \pi$ is defined by $$\ap \pi (x) = \pi (\ap ^{-1}(x))$$ for all $x \in K$.  
It can be easily seen that if $\pi$ and $\pi '$ are equivalent unitary 
representations, then $\ap (\pi )$ and $\ap (\pi ')$ are equivalent.  
Let $[\pi ]$ be denote the equivalent class of unitary representations containing 
$\pi$.   Define  
$\ap [\pi ]=[\ap (\pi )]$ for any automorphism $\ap$ on $K$ and any unitary 
representation $\pi$ of $K$.  Let $\hat K$ be denote the equivalent classes of 
continuous irreducible unitary representations of $K$.  Then $(\ap , [\pi ]) 
\mapsto \ap [\pi ]=[\ap (\pi )]$ defines an action of Aut$~(K)$ on $\hat K$ and 
this action is known as the dual action.  

We now extend Berend's criterion Theorem 2.1 of \cite{Be} to 
finite-dimensional unitary representations.  

\bl\label{ur}
Let $K$ be a compact group and $\Ga$ be a group of automorphisms of $K$.  Then 
the following are equivalent:

\be

\item [(i)] $\Ga$ is ergodic on $K$;

\item [(ii)] $\Ga [\pi ]$ is infinite for any non-trivial $[\pi ]\in \hat K$;

\item [(iii)] $\Ga [\pi ]$ is infinite for any non-trivial finite-dimensional 
unitary representation $\pi$ of $K$.
\ee
\el

\bo
Equivalence of (i) and (ii) is proved in Theorem 2.1 of \cite{Be} and since all 
irreducible unitary representation of compact groups are finite-dimensional, 
(iii) implies (ii).  So, it is sufficient to prove that (ii) implies (iii).  
Let $\pi$ be a finite-dimensional unitary representation of $K$ such that 
$\Ga [\pi ]$ is finite. Let $\Ga _1$ be the subgroup of $\Ga$ consisting of all 
$\ap \in \Ga$ such that $\ap (\pi )$ is equivalent to $\pi$.  
Then $\Ga _1$ is a subgroup of finite index.  

Let $\pi = \oplus _{1\leq j \leq m} \pi _j^{r_j}$ be the decomposition of 
$\pi$ into irreducible unitary representations $\pi _j$ of $K$.  
Then for any automorphism $\ap$ of $K$, $\ap (\pi ) = \oplus _{1\leq j 
\leq m} \ap (\pi _j )^{r_j}$ and each $\ap (\pi _j)$ is also irreducible.  By 
uniqueness of decomposition, each $\ap \in \Ga _1$, defines a permutation on the 
set of equivalent classes of $\pi _j$, $1\leq j \leq m$ 
(cf. 27.30 of \cite{HR2}).  Thus, there is a subgroup $\Ga _2$ of finite 
index in $\Ga _1$ such that $\ap (\pi _j)$ is equivalent to 
$\pi _j$ for all $\ap \in \Ga _2$ and all $j$.  Since $\Ga _1$ has finite 
index in $\Ga$ and $\Ga _2$ has finite index in $\Ga _1$, we get that 
$\Ga _2$ has finite index in $\Ga$.  Using (ii) this implies that 
each $\pi _j$ is trivial and hence $\pi$ is trivial.  Thus, 
(ii) implies (iii) is proved. 
\eo

The following shows stage ergodicity implying ergodicity. 

\bl\label{se}
Let $K$ be a compact group and $\Ga$ be a group of automorphisms of $K$. 
Suppose there is a decreasing sequence $(K_i)_{i\geq 0}$ of closed 
$\Ga$-invariant subgroups with $K_0 =K$ and $\cap K_i = (e)$ such that  
$K_i$ is normal in $K_{i-1}$ and $\Ga$ is ergodic on $K_{i-1}/K_i$.  
Then $\Ga$ is ergodic on $K$. 
\el

\bo
Let $\pi$ be an irreducible unitary representation of $K$ such that 
$\Ga [\pi ]$ is finite.  Then $\pi$ is finite-dimensional.  
Suppose $\pi$ is trivial on $K_i$.  Then let $\tilde \pi $ be the 
restriction of $\pi$ to $K_{i-1}$.  Then $\tilde \pi$ is a unitary representation 
of $K_{i-1}/K_i$ of finite-dimension and $\Ga [\tilde \pi]$ is finite.  Since 
$\Ga$ is ergodic on 
$K_{i-1}/K_i$, it follows from Lemma \ref{ur} that $\tilde \pi$ is trivial.  
Thus, $\pi $ is trivial on $K_{i-1}$ if $\pi$ is trivial on $K_i$ and hence it is 
sufficient to claim that $\pi$ is trivial on some $K_i$.  

Since $\pi$ 
is finite-dimensional, there exists $n\geq 1$ such that $\pi (K)$ is a compact and 
hence a closed subgroup of $GL_n(\C )$.  Then by Cartan's Theorem $\pi (K)$ is a 
real Lie group (cf. Part II, Chapter V, Section 9, \cite{Se}).  Since each $K_i$ 
is compact, we get 
that $(\pi (K_i))$ is a decreasing sequence of closed subgroups of $\pi (K)$.  
Considering dimension and using the fact that compact real Lie groups have only 
finitely many connected components, we get that there exists $i_0$ such that  
$\pi (K_i)= K_{i_0}$ for all $i\geq i_0$.  Since $\cap K_i = (e)$, $\pi (K_i)$ is 
trivial for all $i \geq i_0$.  Thus, $\pi$ is trivial on $K_{i_0}$.
\eo

\bl\label{un}
Let $K$ be a compact group and $\Ga$ be a group of automorphisms of $K$. 
Suppose there exists a collection $(K_i)_{i\in I}$ of closed $\Ga$-invariant 
subgroups $K$ such that $\ol {\cup _{i\in I} K_i }= K$ and $\Ga$ is ergodic on 
$K_i$ for all $i \in I$.  Then $\Ga$ is ergodic on $K$. 
\el

\bo
Let $\pi$ be an irreducible unitary representation of $K$ such that $\Ga [\pi ]$ 
is finite.  For $i \in I$, let $\pi _i$ be the restriction of $\pi$ to $K_i$.  
Then $\pi _i$ is a finite-dimensional unitary representation of $K_i$ and 
$\Ga [\pi _i] $ is finite for all $i \in I$.  Since $\Ga$ is ergodic on $K_i$, 
by Lemma \ref{ur}, we get that $\pi _i$ is trivial and hence $\pi$ is trivial on 
$K_i$ for all $i \in I$.  Since $\ol {\cup K_i }=K$, $\pi$ is trivial on $K$.
\eo

\bl\label{ec}
Let $K$ be a compact group and $\Ga$ be a group of automorphisms of $K$.  
Suppose there exist subgroups $\Ga _1$ and $\Ga _2$ of $\Ga$ such that $\Ga _2$ 
is a compact normal subgroup of $\Ga$ and $\Ga = \Ga _1 \Ga _2$.  Then we have 
the following:

\be

\item [(i)] $\Ga _2[\pi]$ is finite for any irreducible unitary representation 
$\pi$ of $K$; 

\item [(ii)] $\Ga$ is ergodic on $K$ if and only if $\Ga _1$ is ergodic on $K$. 
\ee
\el

\bo
Let $\pi $ be an irreducible unitary representation of $K$.  Let $\cal H$ be the 
Hilbert space on which $\pi$ is defined.  Let $n$ be the dimension of $\cal H$.  
Let $R_n$ be the space of all irreducible unitary 
representations of $K$ of dimension $n$ and $A_n$ be the equivalent classes of 
irreducible unitary representations of $K$ 
of dimension $n$.  Equip $R_n$ with the smallest topology for which the functions 
$\pi \mapsto <\pi (x) u, v>$ are continuous on $R_n$ 
for all $x \in K$ and all $u, v\in {\cal H}$.  
Consider the quotient map $\pi \mapsto [\pi ]$ of $R_n$ 
onto $A_n$ and equip $A_n$ with the quotient topology from $R_n$.  It is 
well-known that $A_n$ is discrete (cf. 3.5.8 and 18.4.3 of \cite{Di}).  
Let $\ap _n \to \ap $ in Aut$~(K)$.  Then for $x \in K$ and $u, v\in {\cal H}$, 
we have 
$$<\ap _n \pi (x)u, v> = <\pi (\ap _n^{-1}x)u ,v> \ra <\pi (\ap ^{-1}x)u, v> = 
<\ap \pi (x)u, v>$$ 
as $\pi$ is continuous and $\ap _n^{-1}(x) \to \ap ^{-1}(x)$.  Thus, 
$\ap _n[\pi ]\to \ap [\pi ]$ in $A_n$.  Since $\Ga _2$ is compact, we get that 
$\Ga _2[\pi ]$ is compact in $A_n$.  Since $A_n$ is discrete, $\Ga _2[\pi ]$ is 
finite.  Thus, proving (i).

Suppose $\Ga $ is ergodic on $K$.  If $\pi$ is an irreducible unitary 
representation of $K$ such that $\Ga _1[\pi ]$ is finite, then since $\Ga _2$ is 
normal, (i) implies that $\Ga [\pi ]= \Ga _2\Ga _1[\pi ]$ is finite.  Since 
$\Ga$ is ergodic, $\pi$ is trivial.  Thus, $\Ga _1$ is ergodic.  This proves (ii).
\eo

We next prove the existence of largest ergodic subgroup for $\Ga$-actions. 

\bp\label{i}
Let $K$ be a compact group and $\Ga$ be a group of automorphisms of $K$.  Then 
there is a $Z(\Ga )$-invariant (largest) closed normal subgroup $L$ of $K$ such 
that $\Ga$ is ergodic on $L$ and $\Ga$ is distal on $K/L$. 
\ep

\bo
Let ${\cal A}$ be the collection of all closed normal $\Ga$-invariant subgroups 
$L$ such that $\Ga$ is ergodic on $L$ and elements of $\cal A$ are ordered by 
inclusion.  Since the trivial subgroup is in $\cal A$, $\cal A$ is non-empty.  
Let $\cal B$ be a chain in $\cal A$.  Let $M= \ol {\cup _{L\in {\cal B}} L }$.  
Since $\cal B$ is a chain, $M$ is a closed normal $\Ga$-invariant 
subgroup of $K$.  By Lemma \ref{un}, $\Ga$ is ergodic on $M$ and hence $M\in 
{\cal A}$ and $L\subset M$ for all $L \in {\cal B}$.  Thus, every chain in 
$\cal A$ has upper bound.  By Zorn's Lemma $\cal A$ has a maximal element and let 
$N$ be the maximal element of $\cal A$.  

Let Inn~$(K)$ be the group of all inner automorphisms of $K$.  Then Inn$~(K)$ is 
a compact normal subgroup of Aut~$(K)$.  Let $\te \Ga $ be the group generated by 
$\Ga$ and all inner automorphisms of $K$.  Then $N$ is $\te \Ga$-invariant and 
$\te \Ga  = \Ga {\rm Inn}~(K)$.  Suppose $\te N$ is a closed 
$\te \Ga $-invariant subgroup of $K$ containing $N$ such that $\te \Ga$ is 
ergodic on $\te N/N$.  Since $\Ga $ is ergodic on $N$ and $\Ga \subset \te \Ga $, 
$\te \Ga$  is ergodic on $N$.  Thus, $\te \Ga $ is ergodic on $\te N$ 
(cf. Lemma \ref{se}).  Now Lemma \ref{ec} implies that $\Ga$ is ergodic on 
$\te N$.  Since $\te N$ is $\te \Ga$-invariant, $\te N$ is normal in $K$ and 
hence $\te N \in {\cal A}$ and $N\subset \te N$.  Since $N$ is a 
maximal element in $\cal A$, $N= \te N$.  By Proposition 2.1 of \cite{Ra}, we get 
that $\te \Ga$ is distal on $K/N$.  Since $\Ga \subset \te \Ga$, $\Ga$ is distal 
on $K/N$.

Let $\ba \in Z(\Ga )$.  Since $\Ga$ is ergodic on $N$, by Theorem 2.1 of 
\cite{Be}, there exists a $x \in N$ such that $\Ga (x)$ is dense in $N$.  
So, there exists a sequence $(\ap _n)$ in $\Ga$ such that $\ap _n(x) \ra e$.  
Now $\ap _n(\ba (x) ) = \ba \ap _n (x) \ra e$.  Since $\Ga$ is distal on 
$K/N$, $\ba (x) \in N$.  Since $N$ is $\Ga$-invariant, this implies that 
$\ba (\ap (x)) = \ap (\ba (x))\in \ap (N)=N$ for all $\ap \in \Ga$.  
Since $\Ga (x)$ is dense in $N$, $\ba (N)\subset N$.  
Thus, $N$ is $Z(\Ga )$-invariant.  

If $H$ is a closed $\Ga$-invariant subgroup of $K$ such that $\Ga$ is ergodic on $H$, 
then since $N$ is normal in $K$, $HN$ is also a closed $\Ga$-invariant subgroup 
and hence $\Ga$ is ergodic on the closed subgroup $HN$ (cf. Lemma \ref{se}).  
This implies that $\Ga$ is ergodic on $HN/N$ but $\Ga$ is distal on $K/N$ and 
hence distal on $HN/N$.  Thus, $N=HN\supset H$.   
\eo

\end{section}

\begin{section}{Distal actions}

In this section we prove that distality of generating elements is sufficient for 
distality of finitely generated abelian group actions on compact groups which is 
needed for our purpose and which may be of independent interest.  

We first consider zero-dimensional compact groups.  It may be noted that 
a compact group is zero-dimensional if and only if it is totally disconnected 
(see 3.5 and 7.7 of \cite{HR1}).  

\bl\label{l1}
Let $\Ga$ be a finitely generated abelian group of automorphisms of a compact 
zero-dimensional group $K$.  Suppose $\ap _1, \ap _2, \cdots , \ap _m \in \Ga$ 
generate $\Ga $.  Then the following are equivalent:

\be
\item [(i)] $\Ga$ action on $K$ is distal;

\item [(ii)] $\ap _i$ is distal on $K$ for $i =1, 2, \cdots , m$; 

\item [(iii)] $K$ has arbitrarily small $\Ga$-invariant compact open subgroups.
\ee
\el

\bo
It is sufficient to prove that (ii) implies (iii).  We prove (iii) by induction 
on the number of generators of $\Ga$.  
If $\Ga$ has only one generator, the result follows from Proposition 
2.1 of \cite{JR} (see also \cite{BW}).  
If $m >1$, then by induction assumption, for any compact open subgroup $W$ of 
$K$, there exists a compact open subgroup $U\subset W$ that is 
invariant under $\ap _1, \ap _2, \cdots ,\ap _{n-1}$.  By Proposition 
2.1 of \cite{JR}, there 
exists a $k$ such that $\cap _{i=0}^k \ap _n^i (U)$ is $\ap _n$-invariant.  
Let $V= \cap _{i=0}^k \ap _n ^i(U)$.  Then $V$ is invariant under $\ap _1, \ap 
_2, \cdots ,\ap _{n-1}$ as $\ap _i$s commute and $V$ is invariant under $\ap 
_n$ also, hence $V$ is invariant under the group generated 
by $\ap _1, \ap _2, \cdots ,\ap _{n-1}$ and $\ap _n$ which is $\Ga$.  Thus, 
proving (iii).
\eo     

We next consider abelian groups and we recall a few results on duality of 
locally compact abelian groups.  Let $G$ be a locally compact abelian group.  
Then any continuous homomorphism of $G$ into the circle group $\{ z\in \C\mid |z| 
=1 \}$ is known as character of $G$.  Let $\hat G$ be the group of characters on 
$G$.  Then $\hat G$ with compact-open topology is known as the dual group of $G$.  
It is well-known that $\hat G$ is discrete if 
and only if $G$ is compact, $\widehat {(\hat G)} \simeq 
G$, and $G$ is metrizable if and only if $\hat G$ is $\sigma$-compact
(cf. Theorems 12, 23 and 29 of \cite{Mo}).  Thus, $K$ is a compact 
(metrizable) abelian group if and only if $\hat K$ is a countable discrete group.    

For any closed subgroup $M$ of $G$, consider the subgroup $M^\perp$ of $\hat G$ 
defined by 
$M^\perp = \{ \chi \in \hat G \mid \chi ~~{\rm is ~~ trivial ~~ on }~~M \}$.   
It follows from Pontryagin duality theory that $\hat M \simeq \hat 
G/M^\perp$ and $M^\perp \simeq \widehat{(G/M)}$ (cf. Theorem 27 of \cite{Mo}).  

For any bi-continuous automorphism $\ap$ of $G$ and $\chi \in \hat G$, 
we define $\ap \chi$ by $\ap \chi (g) = \chi (\ap ^{-1}(g))$ for all $g \in G$ 
and hence $\chi \mapsto \ap \chi$ defines a bi-continuous automorphism of $\hat G$ 
which would be called dual automorphism of $\ap$.  Also, 
$(\ap , \chi ) \mapsto \ap \chi$ defines an action of Aut$~(G)$ on $\hat G$ and it 
is known as the dual action.  

We now apply duality theory of compact abelian groups to obtain the following.

\bl\label{gl}
Let $K$ be a compact abelian group and $\Ga$ be a group of automorphisms 
of $K$.  Suppose there is compact normal subgroup $\Da$ of $\Ga$ and 
$\ap _1, \cdots , \ap _n$ are in $\Ga$ such that $\Ga /\Da$ is abelian and 
$\Ga$ is generated by $\Da$ and $\ap _1, \cdots , \ap _n$.  Then the following 
are equivalent:

\be
\item $\Ga$ is distal on $K$;

\item each $\ap _i$ is distal on $K$.
\ee
\el

\bo
Assume that each $\ap _i$ is distal on $K$.  Proof is based on induction on $n$.  
Suppose $L$ is a non-trivial $\Ga$-invariant subgroup of $K$.  Let $A$ be the 
dual of $L$.  Let $\Ga _1$ be the group generated by $\Da$ and $\ap _1 , \cdots , 
\ap _{n-1}$.  Then by induction hypothesis, $\Ga _1$ is distal on $K$ and hence 
the subgroup $A_1 = \{ \chi \in A \mid \Ga _1(\chi ) ~~{\rm is ~~ finite}~~ \}$ 
is 
non-trivial (cf. Proposition 2.1 of \cite{Ra}).  Since $\Ga /\Da$ is abelian and 
$\Ga _1$ contains $\Da$, we get that $\ap _n$ normalizes $\Ga _1$ and hence $A_1$ 
is $\ap _n$-invariant.  Since $\ap _n$ is distal on $K$, there is a non-trivial 
$\chi  \in A_1$ such that $\ap _n ^m (\chi ) = \chi $ for some $m \not = 0$ 
(cf. Proposition 2.1 of \cite{Ra}).  This implies that 
$\Ga _1 \ap _n^j (\chi ) \subset \cup _{i=1}^{m} \Ga _1(\ap _n^i(\chi))$ 
for all $j$ and hence since $\chi \in A_1$ and $A_1$ is $\ap _n$-invariant, the 
orbit $\Ga (\chi)$ is finite.  Thus, $\Ga$ is not ergodic on $L$.  It 
follows from Proposition 2.1 of \cite{Ra} that $\Ga$ is distal on $K$.
\eo

We now consider connected groups. 

\bl\label{cc}
Suppose $\ap _1, \cdots , \ap _n \in {\rm Aut}~(K)$ are commuting automorphisms 
on a compact connected group $K$.  Then the following are equivalent:

\be
\item  each $\ap _i$ is distal on $K$;

\item the group generated by $\ap _1, \cdots , \ap _n$ is distal on $K$.
\ee
\el

\bo
Assume that each $\ap _i$ is distal on $K$.  
Let $\Ga$ be the group generated by $\ap _1 , \ap _2, \cdots , \ap _n$.  
Let $x \in K$ be such that $e$ is in the closure 
of $\Ga (x)$.  Let $T$ be a maximal pro-torus in $K$ containing $x$ (cf. 
Theorem 9.32 of \cite{HM}).  Then Aut~$(K) = {\rm Inn}~(K)\Omega$ where Inn$~(K)$ 
is the group of inner automorphisms of $K$ and $\Omega = 
\{ \tau \in {\rm Aut}~(K) \mid \tau (T) =T \}$ (cf. Corollary 9.87 of \cite{HM}).  
Let $\ba _i \in \Omega $ be such that $\ap _i \in {\rm Inn}~(K)\ba _i$.  
If $y \in K$ is such that $\ba _i^{k_n} (y) \to e$, then for each $n \geq 1$, 
$\ap _i ^{k_n} = \sigma _{k_n} \ba _i^{k_n} $ for some $\sigma_{k_n} \in 
{\rm Inn}~(K)$ as Inn$~(K)$ is normal in Aut$~(K)$.
Since $K$ is compact, ${\rm Inn}~(K)$ is compact, hence by passing to a 
subsequence, if necessary, we may assume that 
$\sigma _{k_n}\to \sigma \in {\rm Inn}~(K)$ and hence $\ap _i^{k_n} (y) = 
\sigma _{k_n}\ba _i^{k_n}(y) \to \sigma (e)=e$.  
Since $\ap _i$ is distal, $y=e$.  Thus, $\ba _i$ is distal on $K$ and hence in 
particular, $\ba _i$ is distal on $T$.  

Let $\Omega '$ be the group generated $\ba _1, \cdots , \ba _n$.  
Since Inn$~(K)$ is normal, Inn$~(K) \Omega '$ is a group containing $\ap _i$.  
Also $\Omega ' \subset {\rm Inn}~(K)\Ga$, we get that 
${\rm Inn}~(K) \Ga= {\rm Inn}~(K)\Omega '$.  This implies that 
$$\Omega '/ {\rm Inn}~(K)\cap \Omega ' \simeq {\rm Inn}~(K)\Omega '/{\rm Inn}~(K) 
\simeq \Ga {\rm Inn}~(K)/{\rm Inn}~(K)\simeq \Ga /\Ga \cap {\rm Inn}~(K) .$$  
Since $\Ga$ is abelian, $\Omega '/ {\rm Inn}~(K)\cap \Omega'$ is abelian.   
It follows since each $\ba _i$ is distal on $T$ and from Lemma \ref{gl} 
that $\Omega '$ is distal on $T$.

Since $e$ is in the closure of $\Ga (x)$, there exists 
$\eta _n \in \Ga$ such that $\eta _n (x) \to e$.  Since $\Ga \subset
{\rm Inn}~(K)\Omega '$, there exist $\tau _n\in {\rm Inn}~(K)$ and $\ga _n \in 
\Omega '$ such that $\eta _n = \tau _n \ga _n$.  Since ${\rm Inn}~(K)$ is 
compact, we may assume by passing to a subsequence, if necessary that  
$\tau _n \to \tau \in {\rm Inn}~(K)$.  Then $\ga _n(x) = \tau _n ^{-1} 
( \eta _n(x)) \to \tau ^{-1}(e) =e$.  Since $x \in T$, $\ga _n\in \Omega '$ 
and $\Omega '$ is 
distal on $T$, we get that $x =e$.   This shows that $\Ga$ is distal on $K$
\eo

We now consider the general case.  

\bl\label{c1}
Suppose $\ap _1, \cdots , \ap _n \in {\rm Aut}~(K)$ are commuting automorphisms 
on a compact group $K$.  Then the following are equivalent:

\be
\item  each $\ap _i$ is distal on $K$;

\item the group generated by $\ap _1, \cdots , \ap _n$ is distal on $K$.
\ee
\el

\bo
Suppose each $\ap _i$ is distal on $K$.  
Let $\Ga$ be the group generated by $\ap _1 , \ap _2, \cdots , \ap _n$.  
Let $x \in K$ be such that $e$ is in the closure of $\Ga (x)$.  Let $K_0$ be the 
connected component of $e$ in $K$.  Then $K_0$ is a closed normal $\Ga$-invariant 
subgroup of $K$ and $K/K_0$ is a 
zero-dimensional compact group.  Each $\ap _i$ is distal on $K$ implies that each 
$\ap _i$ is distal on $K/K_0$.  It follows from Lemma \ref{l1} that $\Ga$ is 
distal on $K/K_0$.  Since $e$ is in the closure of $\Ga (x)$, $x \in K_0$.  
Now it follows from Lemma \ref{cc} that $x =e$.  
\eo

The following example shows that Lemma \ref{c1} is not true for any acting group.   

\bex 
Let $M$ be any compact group and $\Ga$ be a countably infinite group.  Consider 
the shift 
action of $\Ga$ on $M^\Ga$ defined for $\ap \in \Ga$ and $f\in M^\Ga$ by 
$\ap f(\ba ) = f(\ap ^{-1}\ba )$ for all $\ba \in \Ga$.  Suppose 
$f \in M^\Ga$ is such that $f(\ba )$ is 
identity in $M$ for all but finitely many $\ba \in \Ga$.  We claim that identity 
in $M^\Ga$ is a limit point of $\Ga (f)$.  Let $\ba _1 , \cdots , \ba _n $ be 
such that $f(\ba )$ is identity in $M$ if $\ba \not = \ba _i$.  
Let $U$ be a neighborhood of identity in $M$.  Then for any $\ap _1 , 
\cdots ,\ap _m$, since $\Ga$ is infinite, there is a $\ap \in \Ga$ such that 
$\ap ^{-1}\ap _i \not \in \{ \ba _i \}$ and hence $\ap f(\ap _i)\in U$.  This 
shows that the shift action of $\Ga$ on $M^\Ga$ is not distal.  If we choose 
$\Ga$ to be finitely generated, infinite and torsion, we get that each $\ap \in 
\Ga$ is distal but $\Ga$ is not distal (cf. \cite{Go} for existence of such $\Ga$.)
\eex

\end{section}

\begin{section}{Berend's Theorem}

We now prove Berend's Theorem for ergodic $\Z ^d$-actions on 
compact groups provided that center of the action has DCC.

\bt\label{bt}
Let $K$ be a compact group and $\Ga$ be a finitely generated abelian group of 
automorphisms of $K$.  Suppose $\Ga$ is ergodic on $K$ and $Z(\Ga )$ has DCC on
$K$.  Then there exists $\ap \in \Ga$ such that $\ap$ is ergodic on $K$.  
Furthermore, if $\Ga$ is generated by $\ap _1 ,\cdots ,\ap _n$, 
then there exist positive integers $i_1 , \cdots , i_n$ such that 
$\ap _1 ^{i_1}\cdots \ap _n^{i_n}$ is ergodic on $K$.
\et

We need the following results. 

\bl\label{bl2}
Let $K$ be a compact group and $\ap$, $\ba$ be two commuting automorphisms of 
$K$.  Suppose $\ap$ is ergodic on $K$ and $\ba$ is distal on $K$.  Then 
$\ap ^i\ba ^j$ is ergodic on $K$ for all integers $i$ and $j$ with $i\not =0$.  
\el

\bo
Since $\ap ^i$ is ergodic on $K$ for any $i\not =0$ and $\ba ^j$ is distal on 
$K$ for all $j$, it is sufficient to prove that $\ap \ba$ is ergodic on $K$.  
Applying Proposition 2.1 to $\ap\ba$ on $K$, we get that there exists a closed 
normal $Z(\ap , \ba )\subset Z(\ap \ba)$-invariant subgroup $K_1$ such that $\ap 
\ba$ is ergodic on $K_1$ and $\ap \ba$ is distal on $K/K_1$.  Since $\ap$ and 
$\ba$ commute, $\ap , \ba \in Z(\ap , \ba )$, hence $K_1$ is invariant under 
both $\ap$ and $\ba$.  Since $\ba$ is distal on $K$, $\ba$ is distal on $K/K_1$.  
By Lemma \ref{c1}, the group generated by both $\ap \ba$ and $\ba$ is distal on 
$K/K_1$.  It is easy to see that $\ap$ is in the group generated by $\ap \ba$ and 
$\ba$, hence $\ap$ is also distal on $K/K_1$ but since $\ap$ is ergodic on 
$K$, $\ap$ is ergodic on $K/K_1$ and hence $K = K_1$ 
(cf. Proposition 2.1 of \cite{Ra}).
\eo

\bl\label{bl1}
Let $\ap$, $\ba$ be two commuting automorphisms of $K$ such that $Z(\ap , \ba)$ 
has DCC.  If $\ap $ is ergodic on $K$, then $\ap ^i \ba$ is ergodic for 
all but finitely many $i>0$.
\el

\bo
We first note that since $\ap$ and $\ba$ commute, $\ap , \ba \in Z(\ap , \ba )$.  
Let $K_0 =K$.  We now inductively construct a decreasing sequence 
$(K_i)_{i\geq 1}$ of $Z(\ap, \ba )$-invariant closed subgroups such that 
$K_i$ is normal in $K_{i-1}$, $\ap ^j\ba$ is ergodic on $K_{i-1}/K_i$ for all 
$j \not = i$ and $\ap ^i \ba$ is ergodic on $K_i$.  

Applying Proposition \ref{i} to $\ap\ba$ on $K$, we get that there exists 
$Z(\ap, \ba) \subset Z(\ap \ba)$-invariant closed subgroup $K_1$ 
such that $\ap \ba$ is ergodic on $K_1$ and $\ap \ba$ is distal on $K/K_1$.  
Since $\ap$ is ergodic on $K$, $\ap$ is ergodic on $K/K_1$ and hence by Lemma 
\ref{bl2}, we get that $\ap ^i\ba$ is ergodic on $K/K_1$ for all 
$i \not =1$.  

If for $j\leq i$, $Z(\ap, \ba )$-invariant closed subgroups $K_j$ are chosen so that 
$\ap ^l\ba$ is ergodic on $K_{j-1}/K_j$ for all $l \not =j$ and $\ap ^j \ba$ is 
ergodic on $K_j$.  We now choose $K_{i+1}$.  
Applying Proposition \ref{i} to $\ap ^{i+1}\ba$ on $K_i$, we 
get that there exists a $Z(\ap , \ba )\subset Z(\ap ^{i+1}\ba)$-invariant normal 
subgroup $K_{i+1}$ of 
$K_i$ such that $\ap ^{i+1}\ba$ is ergodic on $K_{i+1}$ and distal on 
$K_i/K_{i+1}$.  Since $\ap ^i\ba$ is ergodic on $K_i$, $\ap ^{i}\ba$ is 
ergodic on $K_i/K_{i+1}$ and hence $\ba ^{-1}\ap ^{-i}$ is ergodic on 
$K_i/K_{i+1}$.   Since $\ap ^{i+1}\ba$ is distal $K_i/K_{i+1}$, 
by Lemma \ref{bl2}, 
$\ap = \ba ^{-1}\ap ^{-i}\ap ^{i+1}\ba$ is ergodic on $K_i/K_{i+1}$.  
Again by Lemma \ref{bl2}, we get that $\ap ^j\ba$ is ergodic on $K_i/K_{i+1}$ for 
all $j\not = i+1$.  Thus, there exists a decreasing sequence $(K_i)_{i\geq 1}$ of 
$Z(\ap, \ba )$-invariant closed subgroups such that 
$K_i$ is normal in $K_{i-1}$, $\ap ^j\ba$ is ergodic on $K_{i-1}/K_i$ for all 
$j \not = i$ and $\ap ^i \ba$ is ergodic on $K_i$.  

Since $Z(\ap, \ba )$ has DCC on $K$, there exists a $i$ such that 
$K_i = K_j$ for all $j > i$.  This implies that for $j >i$, $\ap ^j\ba$ is 
ergodic on $K_{m-1}/K_m$ for $1\leq m\leq i$ and $\ap ^j\ba$ is ergodic on $K_j = 
K_i$.  Thus, by Lemma \ref{se} we get that $\ap ^j\ba$ 
is ergodic on $K$ for all $j>i$.
\eo

\bp\label{b1}
Let $K$ be a compact group and $\Ga$ be a finitely generated group of commuting 
automorphisms of $K$ generated by $\ap_1, \cdots , \ap _n$.  Then there exists 
a series $(e) =K_0 \subset K_{1}\subset \cdots \subset K_n \subset K$ of closed 
normal $Z(\Ga )$-invariant subgroups such that 
\be

\item $\ap _i$ is ergodic on $K_i/K_{i-1}$ and distal on $K/K_i$, 

\item $\Ga$ is ergodic on $K_n$ and distal on $K/K_n$.

\ee
Moreover, if $\Ga$ is ergodic on $K$, $K = K_n$.  
\ep

\bo
We first observe that since $\Ga$ is abelian, $\Ga \subset Z(\Ga )$ and hence any 
$Z(\Ga )$-invariant subgroup is also $\Ga$-invariant.  
It follows from Proposition \ref{i} that there exists a closed normal 
$Z(\Ga )\subset Z(\ap _1)$-invariant subgroup $K_1$ 
such that $\ap _1$ is ergodic on $K_1$ and distal on $K/K_1$.  For $i\geq 1$ and 
for all $1\leq j\leq i$, 
if closed normal $Z(\Ga )$-invariant subgroups $K_j$ are chosen such that 
$\ap _j$ is ergodic on $K_j/K_{j-1}$ and distal on $K/K_j$.  
We now choose 
$K_{i+1}$.  Applying Proposition \ref{i} to $\ap _{i+1}$ on $K/K_i$, there exists 
a closed normal $Z(\Ga )\subset Z(\ap _{i+1})$-invariant subgroup $K_{i+1}$ 
containing $K_i$ such 
that $\ap _{i+1}$ is ergodic on $K_{i+1}/K_i$ and $\ap _{i+1}$ is distal on 
$K/K_{i+1}$.  Thus, we have closed normal $Z(\Ga )$-invariant subgroups 
$(e)=K_0 \subset K_{1}\subset \cdots \subset K_n \subset K$ such that $\ap _i$ is 
ergodic on $K_i/K_{i-1}$ and distal on $K/K_i$.  

Since $\ap _i \in \Ga$, $\Ga$ is ergodic on $K_i/K_{i-1}$ for all $1\leq i \leq 
n$, hence by Lemma \ref{se}, $\Ga$ is ergodic on $K_n$.   

Since $K/K_n$ is a quotient of 
$K/K_i$ for all $1\leq i \leq n$, each $\ap _i$ is distal on $K/K_n$.  Thus by 
Lemma \ref{c1}, $\Ga$ is distal on $K/K_n$.  

If $\Ga$ is ergodic on $K$, then $\Ga$ is ergodic on $K/K_n$ but $\Ga$ is distal 
on $K/K_n$, hence $K=K_n$ (cf. Proposition 2.1 of \cite{Ra}).
\eo

\bo $\!\!\!\!\!$ {\bf of Theorem \ref{bt}}\ \
Let $\ap _1 , \ap _2 , \cdots , \ap _n \in \Ga$ be a set of generators of $\Ga$.  
The proof is by induction on $n$.  Suppose $\Ga$ has more than 
one-generator.  By Proposition \ref{b1}, there exist closed normal 
$\Ga$-invariant subgroups $(e) = K_0 \subset K_1 \subset \cdots \subset K_n = K$ 
such that $\ap _i$ is ergodic on $K_i/K_{i-1}$ and $\ap _i$ is distal on $K/K_i$.

Since $\ap _n$ is ergodic on $K/K_{n-1}$ and $\ap _{n-1}$ is distal on 
$K/K_{n-1}$, Lemma \ref{bl2} implies that $\ap _n^i\ap _{n-1}^j$ is ergodic on 
$K/K_{n-1}$ for all $i$ and $j$ with $i\not = 0$.  In particular, 
$\ap _n \ap _{n-1}^j$ is ergodic on $K/K_{n-1}$ for all $j$. 

Since $Z(\Ga ) \subset Z(\ap _n , \ap _{n-1})$, $Z(\ap _n, \ap _{n-1})$ has DCC 
on $K$.  
Since $\ap _{n-1}$ is ergodic on $K_{n-1}/K_{n-2}$, applying Lemma \ref{bl1} to 
$\ap _{n-1}$ and $\ap _n$ on $K_{n-1}/K_{n-2}$, we get that 
$\ap _n \ap _{n-1}^j$ is ergodic on $K_{n-1}/K_{n-2}$ for all but finitely 
many $j>0$.  Since $\ap _n \ap _{n-1}^j$ is ergodic on $K/K_{n-1}$ for all $j$, 
we get that $\ap _{n-1}^j \ap _n$ is ergodic on $K/K_{n-2}$ 
for some $j>0$ (cf. Lemma \ref{se}).  

Let $\Da$ be the group generated by $\ap _1,\cdots ,\ap _{n-2}, \ap _{n-1}^j
\ap _n$.  Then $Z(\Da ) $ contains $Z(\Ga )$ and hence $Z(\Da )$ has DCC on $K$.  
Since each $\ap _i$ is ergodic on $K_i/K_{i-1}$, we get that $\Da$ is ergodic on 
$K_i/K_{i-1}$ for all $1\leq i \leq n-2$, hence by Lemma \ref{se}, $\Da$ is ergodic 
on $K_{n-2}$.  Since $\ap _{n-1}^j \ap _n $ is ergodic on $K/K_{n-2}$, we get that 
$\Da$ is  ergodic on $K$ (cf. Lemma \ref{se}).  Since $\Da$ is generated by $n-1$ 
elements, result follows from induction assumption. 
\eo

The next result is a consequence of a classical result of Rokhlin \cite{Ro} and 
Theorem \ref{bt}.

\bc\label{rok}
Let $K$ be a compact group and $\Ga$ be a finitely generated group of commuting 
automorphisms of $K$.  Suppose $Z(\Ga)$ has DCC on $K$.  Then the following are 
equivalent:

\be

\item [(1)] $\Ga$ is ergodic on $K$;

\item [(2)] there exists $\ap \in \Ga$ such that $\ap$ is mixing of all orders.

\ee
\ec 

The following results and results in section 5  show that the additional 
condition on $Z(\Ga )$ in Theorem \ref{bt} is not restrictive.

\bc\label{ac1}
Let $K$ be a compact group.  If the center of Aut~$(K)$ has DCC on $K$, then any 
ergodic action of finitely generated abelian group contains ergodic 
automorphisms.
\ec

\bo
Since the center of Aut~$(K)$ is contained in $Z(\Ga )$ for any subgroup $\Ga$ of 
automorphisms, $Z(\Ga )$ has DCC on $K$ if the center of Aut~$(K)$ has DCC on $K$.  
Now the result follows from Theorem \ref{bt}.
\eo

Let $M$ be a compact group 
and $\Ga$ be a countably infinite group.  Consider the shift action of $\Ga$ on 
$M^\Ga$ defined for $\ap \in \Ga $ and $f \in M^\Ga$ by 
$\ap f(\ba) = f(\ap ^{-1}\ba)$ 
for all $\ba \in \Ga$.  In this situation we have the following:

\bt
Let $\Ga$ be a finitely generated abelian group and $M$ be a compact Lie 
group.  Let $K$ be a shift-invariant subgroup $M^\Ga$.  Suppose the action 
of a finitely generated abelian group $\Da$ on $K$ commutes with the shift action 
of $\Ga$ on $K$.  Then $\Da$ is ergodic on $K$ implies that $\Da$ contains 
ergodic automorphisms. 
\et 

\bo
By Theorem 3.2 of \cite{KS} we get that the shift action of $\Ga$ has DCC on 
$M^\Ga$ and hence on $K$.  Since the action of $\Da$ on $K$ commutes with the shift 
action on $K$, $Z(\Da )$ has DCC.  Now the result follows from Theorem \ref{bt}.
\eo

\end{section}

\begin{section}{Some compact abelian groups}

\begin{subsection}{Duals of countable vector spaces}

Let $\F$ be a countable field with discrete topology.  Let $\F ^*$ denote the 
multiplicative group of non-zero elements in $\F$ and $\F ^r$ be the vector space 
over $\F$ of dimension $r\geq 1$.  
Then $\F ^r$ is countable discrete group.  Let $K_r$ be the dual of $\F ^r$.  
Then $K_r$ is a compact group.  Let ${\rm Aut}_{\F }~(\F ^r)$ be the group of all 
$\F$-linear transformation on $\F ^r$ and ${\rm Aut}_{\F}(K_r)$ be the 
corresponding group of dual automorphisms on $K_r$.  In this situation we have 
the following:

\bt\label{f}
Let $\F$ be a countable field with discrete topology.  Let $K _r$ be the dual 
of $\F ^r$ for some integer $r >0$.  Then we have the following: 

\be
\item[(1)] the center of ${\rm Aut}_{\F} ~(K_r)$ has DCC;

\item[(2)] if $\Ga$ is a finitely generated abelian subgroup of 
${\rm Aut}_{\F} ~(K_r)$ such that $\Ga$ is ergodic, then $\Ga$ contains ergodic 
automorphisms.  
\ee
\et

\br
We will now show by an example that finite generation condition on $\Ga$ in 
Theorem \ref{f} (2) can not be relaxed.  Let $p$ be a prime number.  Then for any 
$n \geq 1$, we can find a unique finite field $\F _n$ with $p^n$ number of 
elements 
and $F_n$ contains $F_m$ if and only if $m$ divides $n$ (cf. Corollary 2, Section 1, 
Chapter 1, \cite{We}).  Let $\F = \cup _{k=1}^\infty F_{n^k}$ 
for some fixed $n>1$.  Then $F$ is a countably infinite 
field.  Consider the $\F ^*$-action on $\F$ by scalar multiplication.  Then for 
any non-zero $a \in \F$, $F^*a$ is infinite.  Thus, the dual action of $\F^*$ on 
the dual $K_1$ of $\F$ is ergodic but every element in $\F ^*$ has finite order 
and hence no element of $\F^*$ is ergodic.  
\er

\bo
Consider the $\F ^*$-action on $\F ^r$ given by $q(q_i) = (qq_i)$ for all $q \in 
\F ^*$ and 
$(q_i) \in \F ^r$.  It can be easily seen that $\F ^*$ is in the center of 
${\rm Aut}_{\F }~(\F ^r)$.  If $(L_n)$ is a decreasing sequence of 
closed subgroups in $K_r$ that are $\F ^*$-invariant, then let 
$A_n= \{ \chi \in \F ^r \mid \chi ~~{\rm is ~~ trivial ~~ on }~~ L_n \}$.  Then 
$(A_n)$ is a increasing sequence of subgroups of $\F ^r$ that are $\F ^*$- 
invariant.  Since any subgroup that is $\F ^*$-invariant is a $\F$-subspace of 
$\F ^r$, $(A_n)$ is a increasing sequence of $\F$-subspaces of $\F ^r$.  Since 
$\F ^r$ is finite dimensional, $A_n = A_m$ for all large $n$ and $m$.  
Since $L_n = \widehat {(\F ^r /A_n)}$, we get that $L_n = L_m$ for all large $n$ 
and $m$.  Thus, the center of ${\rm Aut}_{\F }~(K_r)$ has DCC on $K_r$.  

Let $\Ga$ be a finitely generated abelian subgroup of ${\rm Aut}_{\F} ~(K_r)$.  
Then $Z(\Ga )$ contains the center of ${\rm Aut}_{\F} ~(K_r)$.  Since the center 
of ${\rm Aut}_{\F }~(K_r)$ has DCC, we get that $Z(\Ga )$ has DCC.  Now the 
result follows from Theorem \ref{bt}.  
\eo

If we take $K = \Q$, we may obtain Berend's Theorem of \cite{Be} for actions of 
automorphisms as shown below.

\bt\label{rb}
Let $K$ be a compact connected finite-dimensional abelian group and $\Ga$ be a 
group of commuting automorphisms of $K$.  If $\Ga$ is ergodic on $K$, then there 
exists $\ap \in\Ga$ such that $\ap$ is ergodic on $K$.
\et

\bo 
Let $K$ be a compact connected abelian group of finite-dimension.  
Let $r$ be the dimension of $K$.  Let $B_r$ be the dual of the discrete 
group $\Q ^r$.  Then $K$ is a quotient of $B_r$ and any automorphism of $K$ 
lifts to an automorphism of $B_r$.  A group of automorphisms is ergodic on $K$ 
if and only if the corresponding group of lifts is ergodic on 
$B_r$ (cf. \cite{Be} or \cite{Ra}).  Thus, we may assume that $K=B_r$.   

Let $\Ga$ be a group of commuting automorphisms of $B_r$ and $\Ga$ is ergodic on 
$B_r$.  By Lemma 5.9 of \cite{Ra}, there exist a series $(e) =K_0 \subset 
K_1 \subset \cdots \subset K_m = B_r$ of closed connected $\Ga$-invariant 
subgroups with each $K_i \simeq B_{r_i}$ and automorphisms 
$\ap _1 ,\cdots ,\ap _m$ in 
$\Ga$ such that $\ap _i$ is ergodic on $K_i/K_{i-1}$.  This implies that the 
group generated by $\ap _1, \cdots , \ap _m$ is ergodic on $K_m= B_r$ (cf. 
Lemma \ref{se}).  Thus, replacing $\Ga$ by its subgroup generated by 
$\ap _1, \cdots , \ap _m$, if necessary, we may assume 
that $\Ga$ is a finitely generated abelian group.  

It is easy to see that any homomorphism of $\Q ^r$ into $\Q^r$ is a 
$\Q$-linear transformation.  This shows that Aut$~(B_r)\simeq 
{\rm Aut}_{\Q } ~(B_r)$.  Now the result follows from Theorem \ref{f}.
\eo 

\br
Using the results proved here we can prove 
Theorem \ref{rb} for commuting semigroup of epimorphisms of compact connected 
finite-dimensional abelian groups (as lifts of epimorphisms are automorphisms) 
which would retrieve Berend's result in its full generality but as this would cause a 
digression we would not go into the details.  
\er

\end{subsection}

\begin{subsection}{Duals of modules}

We now discuss compact abelian groups arising as duals of modules over rings: 
see \cite{Sh} for any details on rings and modules.  
Let $R$ be a commutative ring and for any $d \geq 1$, let $\bR = \bR _d= 
R [u _1 ^{\pm 1}, \cdots , u_d^{\pm 1}]$ be the 
ring of Laurent polynomials in the commuting variables $u_1, \cdots, u_d$ with 
coefficients in $R$.  Let $\cM$ be a $\bR$-module.  Let $X^{\cM}$ be the dual of 
$\cM$.  Then $X^{\cM}$ is a compact abelian group.  

For $\bn = (n_1, \cdots , n_d) \in \Z ^d$ we define an automorphism $\ap _{\bn}$ 
of $\cM$ by $$\ap _{\bn} a = u_1^{n_1}\dots u_d^{n_d}a$$ for all $a\in \cM$.  
By considering the dual we obtain an automorphism $\ap _{\bn}$ of $X^{\cM}$.  
It can be easily seen that $\bn \mapsto \ap _{\bn}$ defines a $\Z^d$-action on 
$X^{\cM}$.  

Suppose the additive structure of $R$ is cyclic. Then it can be easily seen that 
the afore-defined $\Z ^d$-action has DCC on $X^{\cM}$ if and only if 
$\cM$ is a Noetherian $\bR$-module: a module $A$ over 
a ring $Q$ is called a Noetherian module if any increasing sequence 
$A_1\subset A_2 \subset \cdots \subset A_n\subset A_{n+1}\subset \cdots $ of 
submodules is finite, that is, $A_n =A_m$ for large $n$ and $m$ and a ring $Q$ is 
called a Noetherian ring if $Q$ is a Noetherian $Q$-module.  

If $R= \Z$, Theorem 11.2 of \cite{KS} proved that $\Z ^d$-action defined as above 
has DCC on $X^{\cM}$ if (and only if) $\cM$ is a finitely generated $\bR$-module. 
Theorem 3.3 of \cite{Sc1} has shown that if $R = \Z$, $\cM$ is a finitely 
generated $\bR$-module and the $\Z ^d$-action defined as above is ergodic on 
$X^{\cM}$, then there is a $\bn$ such that $\ap _{\bn}$ is ergodic on $X^{\cM}$.  
In this respect we prove the following using Theorem \ref{bt} 
when the additive group structure of $R$ is a cyclic group.

\bt\label{mo}
Assume that the additive group structure of $R$ is cyclic and $\cM$ is a 
finitely generated $\bR$-module.  Then we have the following:  
\be
\item [(1)] $\{ \ap _{\bn} \mid \bn \in \Z ^d \}$ has DCC on $X^{\cM}$;

\item [(2)] if $\Ga$ is a finitely generated abelian group of $\R$-module 
homomorphisms of $\cM$ and the dual action of $\Ga$ is ergodic on 
$X^{\cM}$, then there is a $\ap \in \Ga$ such that the dual action of 
$\ap$ is ergodic on $X^{\cM}$.
\ee
\et    

\br 
If we take $R=\Z$ and $\Ga = \{ \ap _{\bn} \mid \bn \in \Z ^d \}$, we retrieve 
Theorem 3.3 of \cite{Sc1}.
\er

\bo

We first claim that $\cM$ is Noetherian.  Since the additive structure of $R$ is 
cyclic, $R$ is either finite or $R=\Z$, hence $R$ is Noetherian.  By Hilbert''s
Basis Theorem, the polynomial ring $R[X_1, \cdots , X_d,Y_1, \cdots , Y_d]$ is 
also Noetherian (cf. 8.7 and 8.8 of \cite{Sh}).   
Let $f\colon R[X_1, \cdots , X_d,Y_1, \cdots , Y_d] \to \bR$ be the map defined 
$f(X_i) =u_i$ and $f(Y_i )=u_i^{-1}$ for all $1\leq i \leq d$.  Then $f$ is a 
surjective ring homomorphism 
and hence by Lemma 8.2 of \cite{Sh}, we get that $\bR$ is a Noetherian ring.  
Since $\cM$ is a finitely generated $\bR$-module, we get from Corollary 7.22 (i) 
of \cite{Sh} that $\cM$ is Noetherian.  Since $\Z ^d$-action having DCC on 
$X^{\cM}$ is equivalent to $\cM$ being Noetherian, (1) is proved. 

Let $\ap$ be a $\bR$-module homomorphism.  Then for any $\bn = (n_1, \cdots , 
n_d) \in \Z ^d$, $$\ap \ap _{\bn}(a) = \ap (u_1^{n_1}\cdots u_d^{n_d}a) 
= u_1^{n_1}\cdots u_d^{n_d} \ap (a ) = \ap _\bn \ap (a)$$ 
for all $a \in \cM$.  
Thus, $Z(\Ga )$ contains all $\ap _n$ and hence by (1), 
$Z(\Ga )$ has DCC.  Now (2) follows from Theorem \ref{bt}.
\eo  

\br
Let $\cM$ be free $\bR$-module of rank $k$ (we may take $\cM$ to be direct sum of 
$k$-copies of $\bR$).  Then the group of invertible $\bR$-module homomorphisms of 
$\cM$, can be realized as invertible matrices with coefficients from the ring 
$\bR$.  Thus, $\Ga$ in (2) of Theorem 5.3 may be taken as a finitely generated abelian 
group of invertible matrices with coefficients from $R$.   
\er

We now give one more variation of examples considered in \cite{KS}.  Let $\F$ be 
a countable field and for any $d\geq 1$, let $\bS = \F[u _1^{\pm 1}, \cdots , 
u_d^{\pm 1}]$ be the ring of Laurent polynomials in the commuting variables $u_1, 
\cdots, u_d$ with coefficients in $\F$.  Let $\cN$ be a $\bS$-module.  Let 
$K^{\cN}$ be the dual of $\cN$.  Then $K^{\cN}$ is a compact abelian group.    

Consider the $\Z ^d$-action defined as above, that is, for any 
$\bn = (n_1, \cdots , n_d) \in \Z ^d$ we define an automorphism $\ap _{\bn}$
of $\cN$ by $$\ap _{\bn} a = u_1^{n_1}\dots u_d^{n_d}a$$ for all $a\in \cN$.

Let $\Da = \F ^* \times \Z ^d$.  For $(q , \bn) \in \Da = \F ^* \times \Z^d$, 
we define $\ap _{(q, \bn)}$ by 
$$\ap _{(q, \bn)} (a) = q \ap _{\bn}a $$ for all $a \in \cN$.  By considering the 
dual, we get an automorphism of $K^\cN$.  Thus, we get an action of $\Da$ 
on $K^\cN$.    

It can be easily seen that $\Da$ having DCC on $K^{\cN}$ is equivalent to $\cN$ 
being a Noetherian $\bS$-module.  In this situation we have the following.

\bt\label{mod}
Suppose $\cN$ is finitely generated as a $\bS$-module.  Then we have the following:

\be
\item [(1)] $\Da$ has DCC on $K^{\cN}$, the dual of $\cN$.

\item [(2)] If $\Ga$ is a finitely generated abelian group of $\bS$-module 
homomorphisms of $\cN$ and the dual action of $\Ga$ is ergodic on $K^{\cN}$, 
then there is a $\ap \in \Ga$ such that the dual action of $\ap$ is ergodic on 
$K^{\cN}$.

\item[(3)] If $\Ga = \{ \ap _{\bn} \mid \bn \in \Z ^d \}$, then the dual 
action of $\Ga$ is ergodic on $K^{\cN}$ implies that there is a $\bn \in \Z ^d$ 
such that the dual action of $\ap _{\bn}$ is ergodic on $K^{\cN}$.

\ee
\et

\bo
Since $\Da$ having DCC on $K^{\cN}$ is equivalent to $\cN$ being Noetherian, it 
is sufficient to show that $\cN$ is a Noetherian module over the ring $\bS$ if 
$\cN$ is finitely generated.  

Since $\F$ is a field, it is easy to see that $\F$ is a Noetherian ring as $\F$ 
is the only non-zero ideal.  By Hilbert's Basis Theorem, the polynomial 
ring $F[X_1, \cdots , X_d,Y_1, \cdots , Y_d]$ is also Noetherian (cf. 8.7 and 8.8 
of \cite{Sh}).   Now $X_i \mapsto u_i$ and $Y_i \mapsto u_i^{-1}$ 
induces a ring homomorphism of  $F[X_1, \cdots , X_d,Y_1, \cdots , Y_d]$ onto 
$\bS = \F[u _1^{\pm 1} , \cdots , u_d^{\pm 1}]$ and hence by Lemma 8.2 of 
\cite{Sh}, we get that  $\bS$ is a Noetherian ring.  Since $\cN$ is a finitely 
generated $\bS$-module, Corollary 7.22 (i) of \cite{Sh} implies that $\cN$ is 
Noetherian.  This proves (1).

It is easy to see that any $\bS$-module homomorphism commutes with $\Da$-action 
and hence by (1) $Z(\Ga)$ has DCC if $\Ga$ is a group of $\bS$-module homomorphisms. 
Thus, (2) follows from Theorem \ref{bt}.  

Since each $\ap _{\bn}$ is a $\bS$-module homomorphism, (3) follows from (2).
\eo

We would like to note that compact groups considered in Theorem \ref{f} are 
finite-dimensional whereas compact groups considered in this subsection need not 
be of finite-dimension.  In fact, if $R=\Z$ or $\F = \Q$ and modules are 
free-modules, then compact groups in Theorems \ref{mo} and \ref{mod} are 
infinite-dimensional.

\end{subsection}
\end{section}

\noindent {C. Robinson Edward Raja \newline
Stat-Math Unit \newline
Indian Statistical Instittue \newline
8th Mile Mysore Road \newline
Bangalore 560 059. India \newline
e-mail: creraja@isibang.ac.in}

\end{document}